\documentclass[11pt]{article}
\usepackage{enumerate}
\usepackage{amssymb,a4wide,latexsym,makeidx,epsfig,fleqn}
\usepackage{amsthm}
\usepackage{amsmath}
\usepackage{enumerate}
\newtheorem{theorem}{Theorem}[section]

\newtheorem{lemma}[theorem]{Lemma}

\newtheorem{corollary}[theorem]{Corollary}

\begin{document}
\textwidth 150mm \textheight 225mm
\title{Edge Connectivity, Packing Spanning Trees, and Eigenvalues of Graphs
\thanks{Supported by the National Natural Science Foundation of China (No. 11871398), the Natural Science Basic Research Plan in Shaanxi Province of China (Program No. 2018JM1032) and the Seed Foundation of Innovation and Creation for Graduate Students in Northwestern Polytechnical University (No. ZZ2018171).}}
\author{{Cunxiang Duan, Ligong Wang\footnote{Corresponding author.}, Xiangxiang Liu}\\
{\small Department of Applied Mathematics, School of Science, Northwestern
Polytechnical University,}\\ {\small Xi'an, Shaanxi 710072,
People's Republic
of China.}\\ {\small E-mail: cxduanmath@163.com;lgwangmath@163.com;xxliumath@163.com}\\
}
\date{}
\maketitle
\begin{center}
\begin{minipage}{120mm}
\vskip 0.3cm
\begin{center}
{\small {\bf Abstract}}
\end{center}
{\small Let $\mathcal{G}$ be the set of simple graphs (or multigraphs) $G$ such that for each $G \in \mathcal{G}$  there exists at least two non-empty disjoint proper subsets $V_{1},V_{2}\subseteq V(G)$ satisfying $V(G)\setminus(V_{1} \cup V_{2})\neq \phi$ and edge connectivity $\kappa'(G)=e(V_{i},V(G)\backslash V_{i})$ for $1\leq i \leq 2$. A multigraph is a graph with possible multiple edges, but no loops. Let $\tau(G)$ be the maximum number of edge-disjoint spanning trees of a graph $G$. Motivated by a question of Seymour on the relationship between eigenvalues of a graph $G$ and bounds of $\tau(G)$, we mainly give the relationship between the third largest (signless Laplacian) eigenvalue and the bound of $\kappa'(G)$ and $\tau(G)$ of a simple graph or a multigraph $G\in\mathcal{G}$, respectively.

\vskip 0.1in \noindent {\bf Key Words}: \ eigenvalue, signless Laplacian eigenvalue, edge connectivity, the spanning tree packing number, quotient matrix \vskip
0.1in \noindent {\bf AMS Subject Classification (2016)}: \ 05C50, 05C40, 05C05.}
\end{minipage}
\end{center}

\section{Introduction }
\label{sec:ch6-introduction}

In this paper, we consider finite undirected simple graphs or multigraphs. Undefined notions will follow Bondy and Murty \cite{BoMu}. Let $G$ be a graph with vertex set $V=V(G)=\{v_{1},v_{2},\ldots,v_{n}\}$ and edge set $E=E(G)=\{e_{1},e_{2},\ldots,e_{m}\}$. A  simple graph is a graph without loops and multiple edges. The adjacency matrix of $G$ is an $n$ by $n$ matrix $A(G)=(a_{ij})$, where $a_{ij}$ equals the number of edges between $v_{i}$ and $v_{j}$ for $1\leq i,j \leq n$. If $G$ is simple, then $A(G)$ is a symmetric (0,1)-matrix. Eigenvalues of $A(G)$ are eigenvalues of $G$. We use $\lambda_{i}=\lambda_{i}(G)$ to represent the $i$-th largest eigenvalue of $G$. 
Therefore, we always have $\lambda_{1} \geq \lambda_{2} \geq \cdots \geq \lambda_{n}$. Let $d_{v}$ denote the degree of a vertex $v$ of $G$. Let $\delta$ and $\Delta$ denote the minimum degree and the maximum degree of a graph $G$, respectively. Let $D(G)=diag(d_{v_{1}},d_{v_{2}},\ldots,d_{v_{n}})$ be the diagonal matrix of vertex degrees of $G$. The Laplacian matrix and the signless Laplacian matrix are the matrix $L(G)=D(G)-A(G)$ and $Q(G)=D(G)+A(G), $ respectively. We use $\mu_{i}=\mu_{i}(G)$ and $q_{i}=q_{i}(G)$ to represent the $i$-th largest eigenvalue of $L(G)$ and $Q(G)$, respectively.

For a graph $G$, a path of length $l$ is defined to be an alternating sequence of vertices and edges $u_{1}, e_{1}, u_{2},\ldots , u_{l},e_{l}, u_{l+1},$ where $u_{1},\ldots ,u_{l+1}$ are distinct vertices of $G$, $e_{1},\ldots,e_{l}$ are distinct edges of $G$ and $e_{i}=u_{i}u_{i+1}$ for $i=1,2,\ldots,l.$ If there exists a path between any two vertices of $G$, then $G$ is called connected.
For a graph $G$, the spanning tree packing number, denoted by $\tau(G)$, is the maximum number of edge-disjoint spanning trees of a graph $G$. Let $c(G)$ denote the number of components of $G$, and $\kappa'(G)$ denote the edge connectivity of a graph $G$. We use $\overline{d}(G)$ to denote the average degree of $G$. For $U \subseteq V(G),$ let $ \overline{d}(U)$=$\overline{d}_{G}(U)$ be the average degree of all vertices of $U$ in $G,$ and $\overline{d}(G[U])$ be the average degree of all vertices of the induced subgraph $G[U]$. Let $S$ and $T$ be disjoint subsets of $V(G)$. We denote by $E(S,T)$ the set of edges each of which has one vertex in $S$ and the other vertex in $T$, and let $e(S,T)=\mid E(S,T) \mid$ represents the number of elements in $E(S,T).$

Seymour proposed the problem on predicting $\tau(G)$ by means of the eigenvalues. Cioab\v{a} \cite{Cioa} obtained the condition of the second largest eigenvalue such that a $d$-regular graph is $k$-edge-connected. Cioab\v{a} and Wang \cite{CiWo}, partially answered a question of Seymour, obtained a sufficient eigenvalue condition for the existence of $k$ edge-disjoint spanning trees in a regular graph for $k=2,3.$ Gu et al. \cite{GLLY} proved the relationship between the second largest eigenvalue and edge connectivity and the spanning tree packing number for a simple graph, and proposed a more general conjecture about the spanning tree packing number for a simple graph $G$ with minimum degree $\delta\geq 2k \geq4.$ Their results sharpened theorems of Cioab\v{a} and Wong, they gave a partial solution to the conjecture of Cioab\v{a} and Wong and Seymour's problem. Gu \cite{Gu} mainly proved, for a multigraph $G$ with multiplicity $m$, the algebraic connectivity condition for the existence of $k$-edge-connected and $k$ edge-disjoint spanning trees, respectively. They mainly used the second largest eigenvalue to predict $\kappa'(G)$ and $\tau(G).$ Recently, there are many results concerning the condition of the eigenvalue for the existence of $k$ edge-disjoint spanning trees in a simple graph \cite{GuX,LHL,LHGL,LiSh,Wong} and the third largest eigenvalue in graphs \cite{LiZh,Obou}, respectively.

Let $\mathcal{G}$ be the set of simple graphs (or multigraphs) $G$ such that for each $G \in \mathcal{G}$  there exists at least two non-empty disjoint proper subsets $V_{1},V_{2}\subseteq V$ satisfying $V\setminus(V_{1} \cup V_{2})\neq \phi$ and edge connectivity $\kappa'(G)=e(V_{i},V\backslash V_{i})$ for $1\leq i \leq 2$. A multigraph is a graph with possible multiple edges, but no loops. The multiplicity is the maximum number of edges between any pair of vertices. We are interested in the similar problem on the relationship between the third largest (signless Laplacian) eigenvalue and edge connectivity and the spanning tree packing number of a (multi)graph $G\in \mathcal{G}$, respectively.

In Section 2, some known lemmas and preliminary results are given, including eigenvalue interlacing properties and quotient matrices. In Section 3, we obtain the relationship between the third largest (signless Laplacian) eigenvalue and edge connectivity in a simple graph $G\in\mathcal{G}$. In Section 4, we mainly present, on the basis of the section 3, the relationship between the third largest (signless Laplacian) eigenvalue and the spanning tree packing number in a simple graph $G\in\mathcal{G}$. In Section 5, as corollaries, other eigenvalue conditions on edge connectivity and the spanning tree packing number are presented. In Section 6, we respectively investigate edge connectivity and the spanning tree packing number of a multigraph $G\in\mathcal{G}$ from (signless Laplacian) spectral perspective. And we study other eigenvalue conditions on edge connectivity and the spanning tree packing number, respectively.

\section{Preliminaries}
\label{sec:ch-sufficient}

In this section, we shall give some known results which will be be used in our arguments.

Given a partition $ \pi=\{X_{1}, X_{2}, \ldots, X_{t}\}$ of the set $\{1, 2, \ldots, n\}$ and a matrix $M$ whose rows and columns are labeled with elements in $\{1, 2, \ldots, n\}$, $M$ can be expressed as the following partitioned matrix
$$M=
  \left(\begin{array}{ccc} M_{11} & \cdots & M_{1t}
\\
   \vdots & \vdots & \vdots
 \\
   M_{t1} & \cdots & M_{tt}
  \end{array}\right)
$$
with respect to $\pi$. The quotient matrix $M_{\pi}$ of $M$ with respect to $\pi$ is the $t$ by $t$ matrix $(b_{ij})$ such that each entry $b_{ij}$ is the average row sum of $M_{ij}$.
Suppose that we partition $V(G)$ into $t$ nonempty subsets $V_{1},V_{2},\ldots,V_{t}$. We denote this partition by $\pi$. We use $\overline{d_{i}}$ to represent the average degree of $V_{i}$ for $1\leq i\leq t$. The adjacency quotient matrix $A_{\pi}(G)=A(V_{1},V_{2},\ldots,V_{t})$ of $G$ with respect to $\pi$ is a $t$ by $t$ matrix with entries $b_{ij}$,
$$ b_{ij}=\left\{
\begin{array}{ll}
\overline{d_{i}}-\sum _{1 \leq i\neq j \leq t}\frac{e(V_{i},V_{j})}{\mid V_{i}\mid},& \mbox {if}   ~i = j,
\\
\frac{e(V_{i},V_{j})}{\mid V_{i}\mid},& \mbox {if}   ~i \neq j.
\end{array}
\right.$$
If the partition $\pi$ is not specified, we often use $A_{t}$ to denote the adjacency quotient matrix $A_{\pi}(G)$. Similarly, The signless Laplacian quotient matrix $Q_{\pi}(G)=Q(V_{1},V_{2},\ldots,V_{t})$ of $G$ with respect to $\pi$ is a $t$ by $t$ matrix with entries $c_{ij}$,
$$ c_{ij}=\left\{
\begin{array}{ll}
2\overline{d_{i}}-\sum _{1 \leq i\neq j \leq t}\frac{e(V_{i},V_{j})}{\mid V_{i}\mid},& \mbox {if}   ~i = j,
\\
\frac{e(V_{i},V_{j})}{\mid V_{i}\mid},& \mbox {if}   ~i \neq j.
\end{array}
\right.$$
When it is clear from the context, we often use $Q_{t}$ to denote the signless Laplacian quotient matrix $Q_{\pi}(G)$. While $M_{t}$ is a $t$ by $t$ square real matrix, the following is well known from line algebra \cite{Stra}. $$\lambda_{1}(M_{t})+\lambda_{2}(M_{t})+\cdots\cdots+\lambda_{t}(M_{t})=tr(M_{t}).$$

\noindent\begin{lemma}\label{le:ch-1} \cite{Nash,Tutt}
Let $G$ be a connected graph with $E(G)\neq\varnothing$, and let $k>0$ be an integer. Then $\tau(G)\geq k$ if and only if for any $X\subseteq E(G)$, $\mid X\mid \geq k(c(G-X)-1)$.
\end{lemma}

\noindent\begin{lemma}\label{le:ch-2} (Page 223 in \cite{CVRS})
For any graph $G$, we have $2\delta(G)\leq q_{1}(G)\leq2\Delta(G)$. For a connected graph $G$, equality holds in either place if and only if $G$ is regular.\end{lemma}

By the definition of the quotient matrix $A_{t}$, the sum of all entries in the $ith$ row is $\overline{d_{i}}$. Let $\Delta_{\pi}=\max_{1\leq i \leq t}\{\overline{d_{i}}\}$ and $\delta_{\pi}=\min_{1 \leq i \leq t}\{\overline{d_{i}}\}.$

\noindent\begin{lemma}\label{le:ch-3} (\cite{GLLY})
Let $G$ be a connected graph and $\pi$ be a partition of $V(G).$ Then $\delta_{\pi}\leq\lambda_{1}(A_{\pi}(G))\leq\Delta_{\pi}.$ \end{lemma}

\noindent\begin{lemma}\label{le:ch-4}
Let $G$ be a connected graph and $\pi$ be a partition of $V(G).$ Then $2\delta_{\pi}\leq\ q_{1}(Q_{\pi}(G))\\ \leq 2\Delta_{\pi}.$ \end{lemma}

\noindent {\bf Proof.} Let $G$ be a connected graph and $\pi$ be a partition of $V(G). $ Then $\delta_{\pi}\leq\overline{d}_{i}(D_{\pi}(G))\leq\Delta_{\pi}$ for $1\leq i \leq t$. By Lemma 2.3 and the definition of the signless Laplacian matrix $Q(G)$ of $G$, we have $2\delta_{\pi}\leq\ q_{1}(Q_{\pi}(G))\leq 2\Delta_{\pi}.$\hfill$\blacksquare$

Given two real sequences $\theta_{1} \geq \theta_{2} \geq \cdots \geq  \theta_{n}$ and $\eta_{1} \geq \eta_{2} \geq \cdots \geq \eta_{m}$ with $n>m$, the second sequence is said to interlace the first one if $\theta_{i} \geq \eta_{i} \geq \theta_{n-m+i}$, for $i=1,2,\ldots,m$. When the eigenvalues of a $m$ by $m$ matrix $B$ interlace the eigenvalues of an $n$ by $n$ matrix $A$, it means the non-increasing eigenvalue sequence of $B$ interlaces that of $A$.

\noindent\begin{lemma}\label{le:ch-5} (\cite{BrHa,GoRo,Haem})
Let $G$ be a graph, then the eigenvalues of any quotient matrix of $G$ interlace the eigenvalues of $G$.
\end{lemma}

\noindent\begin{lemma}\label{le:ch-6} (\cite{GLLY})
Let $G$ be a connected simple graph with minimum degree $\delta$ and $U$ be a non-empty proper subset of $V(G)$. If $e(U,V \backslash U)\leq\delta-1$, then $\mid U \mid\geq\delta+1$.
\end{lemma}

\noindent\begin{lemma}\label{le:ch-7} (\cite{Gu})
Let $G$ be a connected multigraph with multiplicity $m$ and minimum degree $\delta$, and $U$ be a non-empty proper subset of $V(G)$. If $e(U,V \backslash U)\leq\delta-1$, then $\mid U \mid\geq\max\{\lceil\frac{\delta+1}{m}\rceil,2\}$.
\end{lemma}

\noindent\begin{lemma}\label{le:ch-8} (\cite{Weyl})
Let $B$ and $C$ be Hermitian matrices of order $n$, and let $1 \leq i,j \leq n$. Then

$(i)$ $\lambda_{i}(B)+\lambda_{j}(C) \leq \lambda_{i+j-n}(B+C)$ if $i+j \geq n+1.$

$(ii)$  $\lambda_{i}(B)+\lambda_{j}(C) \geq \lambda_{i+j-n}(B+C)$ if $i+j \leq n+1.$
\end{lemma}

\noindent\begin{lemma}\label{le:ch-9} (\cite{GLLY})
Let $\delta$, $\Delta$, $\lambda_{2}$, $\mu_{n-1}$ and $q_{2}$ be the minimum degree, maximum degree, second largest eigenvalue, second smallest Laplacian eigenvalue and second largest signless Laplacian eigenvalue of a graph $G$, respectively. Then

$(i)$ $\mu_{n-1}+\lambda_{2} \leq \Delta$.

$ (ii)$ $\delta+\lambda_{2} \leq q_{2}$.
\end{lemma}

\noindent\begin{lemma}\label{le:ch-10} Let $\delta$, $\Delta$, $\lambda_{3}$, $\mu_{n-2}$ and $q_{3}$ be the minimum degree, maximum degree, third largest eigenvalue, third smallest Laplacian eigenvalue and third largest signless Laplacian eigenvalue of a graph $G$, respectively. Then

$(i)$ $\mu_{n-2}+\lambda_{3} \leq \Delta$.

$ (ii)$ $\delta+\lambda_{3} \leq q_{3}$.
\end{lemma}

\noindent {\bf Proof.} Let $A(G)$, $D(G)$, $L(G)$, $Q(G)$ be the adjacency matrix, diagonal degree matrix, Laplacian matrix and signless Laplacian matrix of $G$, respectively.

$(i)$ Since $L(G)=D(G)-A(G)$, we have $D(G)=L(G)+A(G)$. By Lemma 2.8 $(i)$, we obtain $ \lambda_{n-2}(L(G))+\lambda_{3}(A(G)) \leq \lambda_{1}(D(G)).$ Thus $\mu_{n-2}+\lambda_{3} \leq \Delta.$

$(ii)$ Since $Q(G)=D(G)+A(G)$, by Lemma 2.8 $(i)$, we obtain $ \lambda_{n}(D(G))+\lambda_{3}(A(G)) \leq \lambda_{3}(Q(G))$. Thus $ \delta+\lambda_{3} \leq q_{3}.$\hfill$\blacksquare$

\noindent\begin{lemma}\label{le:ch-11}
If $b \geq 3$ and $k \geq 1$, then $\frac{2 (b-1)^{2}}{b(b-2)}k-\frac{2(b-1)}{b(b-2)}<3k-1$.
\end{lemma}

\noindent {\bf Proof.} $\frac{2 (b-1)^{2}}{b(b-2)}k-\frac{2(b-1)}{b(b-2)}<3k-1 \Longleftrightarrow 3k-\frac{2 (b-1)^{2}}{b(b-2)}k>1-\frac{2(b-1)}{b(b-2)} \Longleftrightarrow \frac{b^{2}-2b-2}{b(b-2)}k>\frac{b^{2}-4b+2}{b(b-2)}.$
As $b\geq3$ and $k\geq 1$, it suffices to show that $b^{2}-2b-2>b^{2}-4b+2,$ which is equivalent to $b>2$. It is correct for $b\geq3$. This completes the proof.\hfill$\blacksquare$

\noindent\begin{lemma}\label{le:ch-12}
If $b' \geq 2$ and $k \geq 1$, then $\frac{2b'-1 }{b'-1}k-\frac{2}{b'-1}<3k-1$.
\end{lemma}

\noindent {\bf Proof.} $\frac{2b'-1 }{b'-1}k-\frac{2}{b'-1}<3k-1 \Longleftrightarrow 3k-\frac{2b'-1 }{b'-1}k>1-\frac{2}{b'-1} \Longleftrightarrow \frac{b'-2}{b'-1}k>\frac{b'-3}{b'-1}.$
It is obviously correct when $b'\geq2$ and $k\geq 1.$ This completes the proof.\hfill$\blacksquare$

\section{Eigenvalues and edge connectivity}
\label{sec:ch-inco}

In this section, we mainly study the relationship between the third largest (signless Laplacian) eigenvalue and edge connectivity in a simple graph $G\in\mathcal{G}$.

The following result is the relationship between the third largest eigenvalue and edge connectivity in a simple graph $G\in\mathcal{G}$. We prove this result by using the adjacency quotient matrix of a graph.

\noindent\begin{theorem}\label{th:ch-1}
Let $k$ be an integer with $k \geq 2$ and $G\in\mathcal{G}$ be a simple graph with minimum degree $\delta\geq 2k-1.$ If $\lambda_{3}(G)<2\delta-\Delta-\frac{4(k-1)}{\delta+1}$, then $\kappa'(G)\geq k.$
\end{theorem}

\noindent {\bf Proof.} We prove this theorem by contradiction and assume that $\kappa'(G)\leq k-1$. Since $G\in \mathcal{G}$, there exist two non-empty disjoint proper subsets $V_{1},V_{2}\subseteq V(G)$ such that $r_{i}=e(V_{i},V\backslash V_{i})=\kappa'(G)\leq k-1 $ for $1\leq i \leq 2$. Let $V'=V \backslash (V_{1}\cup V_{2})$, then $r'=e(V',V_{1})+e(V',V_{2})\leq 2(k-1) \leq \delta-1.$ By Lemma 2.6, we have $\mid V_{i}\mid\geq\delta+1$ for $1\leq i \leq 2$ and $\mid V'\mid\geq\delta+1$. Let $y=e(V_{1},V_{2})\geq 0.$ The adjacency quotient matrix of $G$ with respect to the partition $(V_{1},V_{2},V')$ is

$$A_{3}=
  \left(\begin{array}{ccc}\overline{d_{1}}-\frac{r_{1}}{\mid V_{1}\mid} & \frac{y}{\mid V_{1}\mid} & \frac{r_{1}-y}{\mid V_{1}\mid}
\\
   \frac{y}{\mid V_{2}\mid} & \overline{d_{2}}-\frac{r_{2}}{\mid V_{2}\mid} & \frac{r_{2}-y}{\mid V_{2}\mid} \\
\frac{r_{1}-y}{\mid V'\mid}&\frac{r_{2}-y}{\mid V'\mid}& \overline{d'}-\frac{r_{1}+r_{2}-2y}{\mid V'\mid}
  \end{array}\right),
$$
where $\overline{d'}$ denotes the average degree of $V'$ in $G$, and $\overline{d_{i}}$ denotes the average degree of $V_{i}$ in $G$ for $i=1,2$.
By Lemma 2.3, we have $\lambda_{1}(A_{3}) \leq \max\{\overline{d_{1}},\overline{d_{2}},\overline{d'}\}.$
By $tr(A_{3})=\lambda_{1}(A_{3})+\lambda_{2}(A_{3})+\lambda_{3}(A_{3})$, we have
\begin{align*}
\lambda_{2}(A_{3})+\lambda_{3}(A_{3})&=tr(A_{3})-\lambda_{1}(A_{3})\geq\overline{d_{1}}+
\overline{d_{2}}+\overline{d'}-\frac{r_{1}}{\mid V_{1}\mid}-\frac{r_{2}}{\mid V_{2}\mid}-
\frac{r_{1}+r_{2}-2y}{\mid V'\mid}\\& - \max\{\overline{d_{1}},\overline{d_{2}},\overline{d'}\}
\geq 2\delta -
\frac{2(r_{1}+r_{2})-2y}{\delta+1} \geq 2\delta-\frac{4(k-1)}{\delta+1}.
\end{align*}
By Lemma 2.5, we have $\lambda_{3}(A(G)) \geq \lambda_{3}(A_{3}).$
Therefore, we have
\begin{align*}
\lambda_{3}(A(G))\geq\lambda_{3}(A_{3})\geq 2\delta-\frac{4(k-1)}{\delta+1}-
\lambda_{2}(A_{3})\geq 2\delta-\frac{4(k-1)}{\delta+1}-\Delta,
\end{align*}
which is contrary to the assumption in Theorem 3.1 that $\lambda_{3}(A(G))<2\delta-\Delta-\frac{4(k-1)}{\delta+1}.$ This completes the proof of the theorem.\hfill$\blacksquare$

By Theorem 3.1, we have the following the corollary.

\noindent\begin{corollary}\label{co:ch-2}
Let $k$ and $d$ be two integers with $d\geq 2k-1 \geq 3$ and $G\in\mathcal{G}$ be a d-regular simple graph. If $\lambda_{3}(G)<d-\frac{4(k-1)}{d+1}$, then $\kappa'(G)\geq k.$
\end{corollary}

The following result is the relationship between the third largest signless Laplacian eigenvalue and edge connectivity in $G\in\mathcal{G}$. We prove this result by using the signless Lapacian quotient matrix of a graph.

\noindent\begin{theorem}\label{th:ch-3}
Let $k$ be an integer with $k \geq 2$ and $G\in\mathcal{G}$ be a simple graph with minimum degree $\delta\geq 2k-1.$ If $q_{3}(G)<4\delta-2\Delta-\frac{4(k-1)}{\delta+1}$, then $\kappa'(G)\geq k.$
\end{theorem}

\noindent {\bf Proof.} We prove this theorem by contradiction and assume that $\kappa'(G)\leq k-1$. Since $G\in \mathcal{G}$, there exist two non-empty disjoint proper subsets $V_{1},V_{2}\subseteq V(G)$ such that $r_{i}=e(V_{i},V\backslash V_{i})=\kappa'(G)\leq k-1 $ for $1\leq i \leq 2$. Let $V'=V \backslash (V_{1}\cup V_{2})$, then $r'=e(V',V_{1})+e(V',V_{2})\leq 2(k-1) \leq \delta-1.$ By Lemma 2.6, we have $\mid V_{i}\mid\geq\delta+1$ for $1\leq i \leq 2$ and $\mid V'\mid\geq\delta+1$. Let $y=e(V_{1},V_{2})\geq 0.$ The signless Laplacian quotient matrix of $G$ with respect to the partition $(V_{1},V_{2},V')$ is

$$Q_{3}=
  \left(\begin{array}{ccc}2\overline{d_{1}}-\frac{r_{1}}{\mid V_{1}\mid} & \frac{y}{\mid V_{1}\mid} & \frac{r_{1}-y}{\mid V_{1}\mid}
\\
   \frac{y}{\mid V_{2}\mid} & 2\overline{d_{2}}-\frac{r_{2}}{\mid V_{2}\mid} & \frac{r_{2}-y}{\mid V_{2}\mid} \\
\frac{r_{1}-y}{\mid V'\mid}&\frac{r_{2}-y}{\mid V'\mid}& 2\overline{d'}-\frac{r_{1}+r_{2}-2y}{\mid V'\mid}
  \end{array}\right),
$$
where $\overline{d'}$ denotes the average degree of $V'$ in $G$, and $\overline{d_{i}}$ denotes the average degree of $V_{i}$ in $G$ for $i=1,2$.
By Lemma 2.4, we have $ q_{1}(Q_{3}) \leq 2\max\{\overline{d_{1}},\overline{d_{2}},\overline{d'}\}.$
By $tr(Q_{3})=q_{1}(Q_{3})+q_{2}(Q_{3})+q_{3}(Q_{3})$, we have
\begin{align*}
q_{2}(Q_{3})+q_{3}(Q_{3})&=tr(Q_{3})-q_{1}(Q_{3})=2\overline{d_{1}}+
2\overline{d_{2}}+ 2\overline{d'}-\frac{r_{1}}{\mid V_{1}\mid}-\frac{r_{2}}{\mid V_{2}\mid}-
\frac{r_{1}+r_{2}-2y}{\mid V'\mid}-q_{1}(Q_{3})\\& \geq 2\overline{d_{1}}+2\overline{d_{2}}+ 2\overline{d'} -
\frac{2(r_{1}+r_{2})-2y}{\delta+1}-2\max\{\overline{d_{1}},\overline{d_{2}},\overline{d'}\} \geq 4\delta-\frac{4(k-1)}{\delta+1}.
\end{align*}
By Lemma 2.5, we have $q_{3}(Q(G)) \geq q_{3}(Q_{3}).$
Therefore, we have
\begin{align*}
q_{3}(Q(G))\geq q_{3}(Q_{3})&\geq 4\delta-\frac{4(k-1)}{\delta+1}-q_{2}(Q_{3})\geq 4\delta-\frac{4(k-1)}{\delta+1}-2\max\{\overline{d_{1}},\overline{d_{2}},\overline{d'}\}\\&\geq 4\delta-\frac{4(k-1)}{\delta+1}-2\Delta,
\end{align*}
which is contrary to the assumption in Theorem 3.3 that $q_{3}(Q(G))<4\delta-2\Delta-\frac{4(k-1)}{\delta+1}.$ This completes the proof of the theorem.\hfill$\blacksquare$

By Theorem 3.3, we have the following corollary.

\noindent\begin{corollary}\label{co:ch-4}
Let $k$ and $d$ be two integers with $d\geq 2k-1 \geq 3$ and $G\in\mathcal{G}$ be a d-regular simple graph. If $q_{3}(G)<2d-\frac{4(k-1)}{d+1}$, then $\kappa'(G)\geq k.$
\end{corollary}

{\bf Remark 1.} Let $G\in\mathcal{G}$ be a graph shown in Figure 1. By Matlab calculation and Theorems 3.1 and 3.3, we have $\lambda_{3}\approx 1.732<2\times3-3-\frac{4(2-1)}{3+1}=2,$ $ q_{3}\approx 4.732<4\times3-2\times3-\frac{4(2-1)}{3+1}=5,$ $\kappa'(G)=3>k=2.$

\begin{center}
\includegraphics [width=4 cm, height=3 cm]{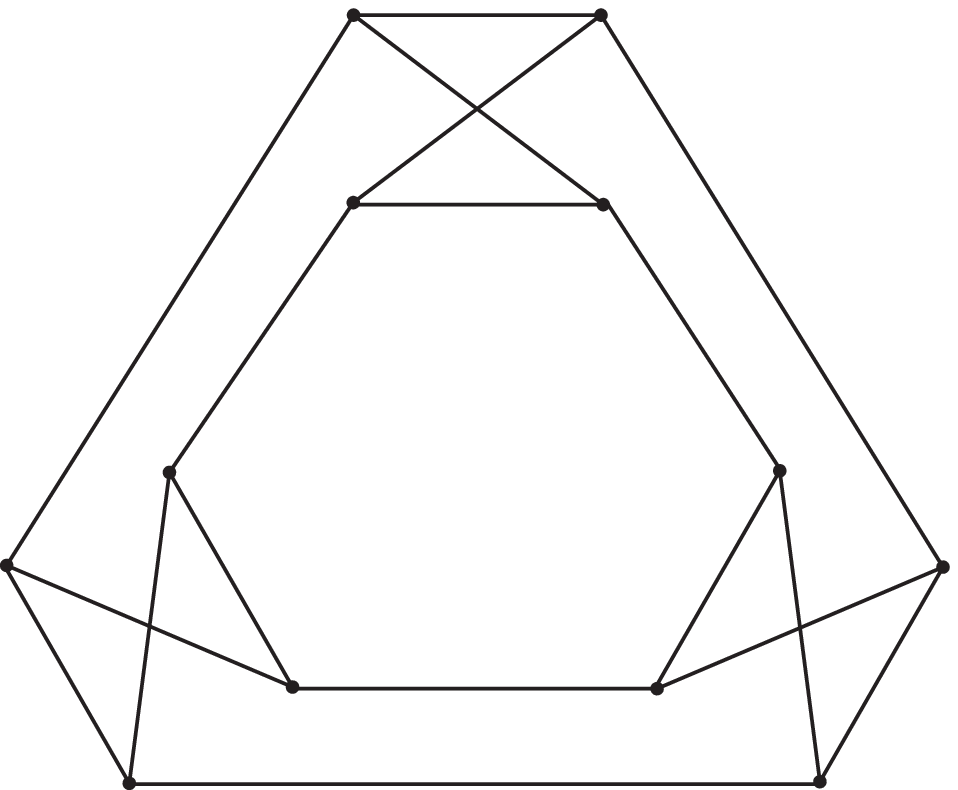}
\centerline{Figure 1: A 3-regular graph with $\lambda_{3}\approx 1.732<2.$}
\end{center}

The following Theorem 3.5 has been given by Gu et al. \cite{GLLY}. In the following, we give a different proof by using the signless Laplacian quotient matrix of a graph.

\noindent\begin{theorem}\label{th:ch-5} (\cite{GLLY})
Let $k$ be an integer with $k \geq 2$ and $G$ be a simple graph with minimum degree $\delta\geq k.$ If $q_{2}(G)<2\delta-\frac{2(k-1)}{\delta+1}$, then $\kappa'(G)\geq k.$
\end{theorem}

\noindent {\bf Proof.} We prove the theorem by contradiction and assume that $\kappa'(G)\leq k-1$. Then there exists a non-empty proper subset $V_{1}\subseteq V(G)$ such that $r=e(V_{1},V\backslash V_{1})\leq k-1 $. Let $V_{2}=V\backslash V_{1}$. By Lemma 2.6, then we have $\mid V_{1}\mid\geq\delta+1$ and $\mid V_{2}\mid\geq\delta+1$. The signless Laplacian quotient matrix of $G$ with respect to the partition $(V_{1},V_{2})$ is

$$Q_{2}=
  \left(\begin{array}{cc}2\overline{d_{1}}-\frac{r}{\mid V_{1}\mid} & \frac{r}{\mid V_{1}\mid}
\\
   \frac{r}{\mid V_{2}\mid} & 2\overline{d_{2}}-\frac{r}{\mid V_{2}\mid}
  \end{array}\right),
$$
where $\overline{d_{i}}$ denotes the average degree of $V_{i}$ in $G$ for $i=1,2$.
By Lemma 2.4, we have $q_{1}(Q_{3}) \leq 2\max\{\overline{d_{1}},\overline{d_{2}}\}.$
By $tr(Q_{2})=q_{1}(Q_{2})+q_{2}(Q_{2})$, we have
\begin{align*}
q_{1}(Q_{2})+q_{2}(Q_{2})&=tr(Q_{2})= 2\overline{d_{1}}+
2\overline{d_{2}}-\frac{r}{\mid V_{1}\mid}-\frac{r}{\mid V_{2}\mid}
 \geq 2\overline{d_{1}}+2\overline{d_{2}} -
\frac{2(k-1)}{\delta+1}.
\end{align*}
By Lemma 2.5, we have $q_{2}(Q(G)) \geq q_{2}(Q_{2}).$
Therefore, we have
\begin{align*}
q_{2}(Q(G))& \geq q_{2}(Q_{2}) \geq 2\overline{d_{1}}+2\overline{d_{2}} -\frac{2(k-1)}{\delta+1}-q_{1}(Q_{2})  \geq 2\overline{d_{1}}+2\overline{d_{2}}-\frac{2(k-1)}{\delta+1}\\&- 2 \max\{\overline{d_{1}},\overline{d_{2}}\} \geq 2\delta-\frac{2(k-1)}{\delta+1},
\end{align*}
which is contrary to the assumption in Theorem 3.5 that $q_{2}(Q(G))<2\delta-\frac{2(k-1)}{\delta+1}.$ This completes the proof of the theorem.\hfill$\blacksquare$

By Theorem 3.5, we have the following the corollary.

\noindent\begin{corollary}\label{co:ch-6}
Let $k$ and $d$ be two integers with $d\geq k \geq 2$ and $G$ be a d-regular simple graph. If $q_{2}(G)<2d-\frac{2(k-1)}{d+1}$, then $\kappa'(G)\geq k.$
\end{corollary}

\section{Eigenvalues and the spanning tree packing number}
\label{sec:ch-fudi}

In this section, we mainly investigate the relationship between the third largest (signless Laplacian) eigenvalue and the spanning tree packing number in a simple graph $G\in\mathcal{G}$.

 By Lemma 2.1, if $\tau(G)\leq k-1,$ then there exists an edge subset $X\subseteq E(G)$ such that $\mid X \mid \leq k(c(G-X)-1)-1$. Let $c(G-X)=s$ and $G_{1},G_{2},\ldots,G_{s}$ be the components of $G-X$. For $1 \leq i \leq s$, let $V_{i}=V(G_{i}),$ $E_{i}=E(G_{i})$, and $r_{i}=e(V_{i},V \backslash V_{i})$. Without loss of generality, we always assume that $$r_{1} \leq r_{2} \leq \cdots \leq r_{s}.\eqno{(1)}$$  With these notations and by $\mid X \mid \leq k(c(G-X)-1)-1$, we have $$\sum _{1 \leq i<j \leq s}e(V_{i},V_{j}) \leq k(s-1)-1=ks-k-1. \eqno{(2)}$$

\noindent\begin{theorem}\label{th:ch-1} For $G\in\mathcal{G}$ and $\delta \geq 2k\geq4$, if $\lambda_{3}(G)< 2\delta-\Delta-\frac{2(3k-1)}{\delta+1}$. Then $r_{i} \geq k$ for any $i$ with $1 \leq i \leq s$.
\end{theorem}

\noindent {\bf Proof.} We argue by contradiction and assume that $r_{i} <k$ for some $i$.   Then $ \kappa'(G) < k$. By Theorem 3.1,
$$\lambda_{3}(G) \geq 2\delta-\Delta-\frac{4(k-1)}{\delta+1} > 2\delta-\Delta-\frac{2(3k-1)}{\delta+1},$$ which is contrary to the assumption that $\lambda_{3}(G)< 2\delta-\Delta-\frac{2(3k-1)}{\delta+1}$.
Therefore, we have $r_{i} \geq k$.\hfill$\blacksquare$

\noindent\begin{theorem}\label{th:ch-2} Let $k$ be an integer with  $k \geq 2$ and $G\in\mathcal{G}$ be a simple graph with minimum degree $\delta \geq 2k $. If $\lambda_{3}(G)< 2\delta-\Delta-\frac{2(3k-1)}{\delta+1}$, then $\tau(G)\geq k$.
\end{theorem}

\noindent {\bf Proof.} We shall argue by contradiction and assume that $\tau(G)\leq k-1.$  By (2) with $k \geq 2$, we have $\sum\limits^{s}_{i=1}r_{i} =2 \sum _{1 \leq i<j \leq s}e(V_{i},V_{j}) \leq 2ks-2(k+1).$ Let $x_{l}$ denote the multiplicity of $l$ in $\{r_{1},r_{2},\ldots,r_{s}\}$ for $1 \leq l \leq 2k-1$. By Theorem 4.1, we have $r_{s} \geq \cdots \geq r_{2} \geq r_{1} \geq k$. Thus $x_{j}=0$ for $j=1,2,\ldots,k-1$. By (2), we have
$$kx_{k}+(k+1)x_{k+1}+\cdots+(2k-1)x_{2k-1}+2k(s-(x_{k}+x_{k+1}+\cdots+x_{2k-1})) \leq \sum\limits^{s}_{i=1}r_{i} \leq 2ks-2(k+1),$$
which implies that
$kx_{k}+(k-1)x_{k+1}+\cdots+2x_{2k-2}+x_{2k-1} \geq 2(k+1).$
Let $h$ be the smallest index such that $x_{h} \neq 0$, then we have
$$ (2k-h)x_{h}+(2k-h-1)x_{h+1}+\cdots+2x_{2k-2}+x_{2k-1} \geq 2(k+1).\eqno{(3)}$$

{\bf Case} 1. $x_{h} \geq 2$.

Since $ h \geq k$, we have $ 2(k+1)>2(2k-h)$. Then there exists an integer $b \geq 3$ such that $(b-1)(2k-h)<2(k+1) \leq b(2k-h)$. By $2(k+1) \leq b(2k-h)$, we have $h \leq \frac{(2b-2)k-2}{b}.$ It follows by $(b-1)(2k-h)<2(k+1)$ and (3) that $x_{h}+x_{h+1}+\cdots+x_{2k-2}+x_{2k-1} \geq b,$ and so by (1), we have $2k-1 \geq r_{b} \geq \cdots \geq r_{2} \geq r_{1}.$ By Lemma 2.6, we have $\mid V_{i} \mid \geq \delta+1$ with $1 \leq i \leq b.$ Let $V'=V \backslash (V_{1} \cup V_{2})$, then $\mid V' \mid \geq \mid V_{3}\mid + \cdots+ \mid V_{b} \mid \geq (b-2)(\delta+1).$ Let $y=e(V_{1},V_{2})\geq 0.$ The adjacency quotient matrix of $G$ with respect to the partition $(V_{p},V_{q},V')$ is
$$A_{3}=
  \left(\begin{array}{ccc}\overline{d_{1}}-\frac{h}{\mid V_{1}\mid} & \frac{y}{\mid V_{1}\mid} & \frac{h-y}{\mid V_{1}\mid}
\\
   \frac{y}{\mid V_{2}\mid} & \overline{d_{2}}-\frac{h}{\mid V_{2}\mid} & \frac{h-y}{\mid V_{2}\mid} \\
\frac{h-y}{\mid V'\mid}&\frac{h-y}{\mid V'\mid}& \overline{d'}-\frac{2(h-y)}{\mid V'\mid}
  \end{array}\right),
$$
where $\overline{d'}$ denotes the average degree of $V'$ in $G$, and $\overline{d_{i}}$ denotes the average degree of $V_{i}$ in $G$ for $i=1,2$.
By Lemma 2.3, we have
$\lambda_{1}(A_{3}) \leq \max\{\overline{d_{1}},\overline{d_{2}},\overline{d'}\}.$ By $tr(A_{3})=\lambda_{1}(A_{3})+\lambda_{2}(A_{3})+\lambda_{3}(A_{3})$ and Lemma 2.11, we have
\begin{align*}
\lambda_{2}(A_{3})+\lambda_{3}(A_{3})&=tr(A_{3})-\lambda_{1}(A_{3})\geq\overline{d_{1}}+
\overline{d_{2}}+\overline{d'}-\frac{h}{\mid V_{1}\mid}-\frac{h}{\mid V_{2}\mid}-
\frac{2(h-y)}{\mid V'\mid} - \max\{\overline{d_{1}},\overline{d_{2}},\overline{d'}\}\\&\geq 2\delta -\frac{2\frac{b-1}{b-2}h}{\delta+1} \geq 2(\delta-\frac{\frac{2 (b-1)^{2}}{b(b-2)}k-\frac{2(b-1)}{b(b-2)}}{\delta+1})>2(\delta-\frac{3k-1}{\delta+1}).
\end{align*}
By Lemma 2.5, we have $\lambda_{3}(A(G)) \geq \lambda_{3}(A_{3})$.
Therefore, we have
$$\lambda_{3}(A(G))\geq\lambda_{3}(A_{3})\geq 2(\delta-\frac{3k-1}{\delta+1})-
\lambda_{2}(A_{3}) \geq 2\delta-\frac{2(3k-1)}{\delta+1}-\Delta,$$
which is contrary to the assumption in Theorem 4.2 that $\lambda_{3}(A(G))<2\delta-\frac{2(3k-1)}{\delta+1}-\Delta$.

{\bf Case} 2. $x_{h} = 1$.

In this case, we can rewrite (3) as $(2k-h-1)x_{h+1}+\cdots+2x_{2k-2}+x_{2k-1} \geq 2(k+1)-(2k-h)=h+2 \geq k+2.$ Let $h'$ be the smallest index such that $x_{h'}>0$ with $h'>h$, then
$$(2k-h')x_{h'}+(2k-h'-1)x_{h'+1}+\cdots+2x_{2k-2}+x_{2k-1} \geq h+2 \geq k+2.\eqno{(4)}$$
Since $ h' > h \geq k$, we have $h'+2>k$. So $k+2>2k-h'$. Thus there exists an integer $b' \geq 2$ such that $(b'-1)(2k-h')<k+2 \leq b'(2k-h')$. By $k+2 \leq b'(2k-h')$, we have $ h'\leq \frac{(2b'-1)k-2}{b'}$. It follows by $(b'-1)(2k-h')<k+2$ and (4) that $x_{h'}+x_{h'+1}+\cdots+x_{2k-2}+x_{2k-1} \geq b',$ and so by (1), we have $2k-1 \geq r_{b'+1} \geq \cdots \geq r_{2} \geq r_{1}. $ By Lemma 2.6, $\mid V_{i} \mid \geq \delta+1$ with $1 \leq i \leq b'+1.$ Let $V'=V \backslash (V_{1} \cup V_{2})$, then $\mid V' \mid \geq \mid V_{3}\mid + \cdots+ \mid V_{b'+1} \mid \geq (b'-1)(\delta+1). $ Let $y=e(V_{1},V_{2})\geq0.$ The adjacency quotient matrix of $G$ with respect to the partition $(V_{p},V_{q},V')$ is
$$A_{3}=
  \left(\begin{array}{ccc}\overline{d_{1}}-\frac{h}{\mid V_{1}\mid} & \frac{y}{\mid V_{1}\mid} & \frac{h-y}{\mid V_{1}\mid}
\\
   \frac{y}{\mid V_{2}\mid} & \overline{d_{2}}-\frac{h'}{\mid V_{2}\mid} & \frac{h'-y}{\mid V_{2}\mid} \\
\frac{h-y}{\mid V'\mid}&\frac{h'-y}{\mid V'\mid}& \overline{d'}-\frac{h'+h-2y}{\mid V'\mid}
  \end{array}\right),
$$
where $\overline{d'}$ denotes the average degree of $V'$ in $G$, and $\overline{d_{i}}$ denotes the average degree of $V_{i}$ in $G$ for $i=1,2$.
By Lemma 2.3, we have $\lambda_{1}(A_{3}) \leq \max\{\overline{d_{1}},\overline{d_{2}},\overline{d'}\}.$
By $tr(A_{3})=\lambda_{1}(A_{3})+\lambda_{2}(A_{3})+\lambda_{3}(A_{3})$ and Lemma 2.12, we have
\begin{align*}
\lambda_{2}(A_{3})+\lambda_{3}(A_{3})&=tr(A_{3})-\lambda_{1}(A_{3})\geq\overline{d_{1}}+
\overline{d_{2}}+\overline{d'}-\frac{h}{\mid V_{1}\mid}-
\frac{h'}{\mid V_{2}\mid}-\frac{h'+h-2y}{\mid V'\mid}\\& - \max\{\overline{d_{1}},\overline{d_{2}},\overline{d'}\}
\geq 2\delta -\frac{b'h+b'h'}{(\delta+1)(b'-1)} \geq 2\delta-\frac{2b'h' }{(\delta+1)(b'-1)}\\&\geq 2(\delta-\frac{\frac{2b'-1}{b'-1}k-\frac{2}{b'-1}}{\delta+1}) >2(\delta-\frac{3k-1}{\delta+1}).
\end{align*}
By Lemma 2.5, we have $\lambda_{3}(A(G)) \geq \lambda_{3}(A_{3})$. Therefore, we have
$$\lambda_{3}(A(G))\geq\lambda_{3}(A_{3})\geq 2(\delta-\frac{3k-1}{\delta+1})-
\lambda_{2}(A_{3}) \geq 2\delta-
\frac{2(3k-1)}{\delta+1}-\Delta,$$
which is contrary to the assumption in Theorem 4.2 that $\lambda_{3}(A(G))<2\delta-\Delta-\frac{2(3k-1)}{\delta+1}$.\hfill$\blacksquare$

\noindent\begin{corollary}\label{th:ch-3} For $k \geq 2, $ and $d \geq 2k$. If $G\in\mathcal{G}$ is a d-regular simple graph with $\lambda_{3}(G)<d-\frac{2(3k-1)}{d+1}$, then $\tau(G)\geq k$.
\end{corollary}

Next, we will discuss the relationship between $q_{3}(G)$ and $\tau(G)$ of a simple graph $G\in\mathcal{G}$.

\noindent\begin{corollary}\label{th:ch-4} For $G\in\mathcal{G}$ and $\delta \geq 2k \geq4$, if $q_{3}(G)< 4\delta-2\Delta-\frac{2(3k-1)}{\delta+1}$. Then $r_{i} \geq k$ for any $i$ with $1 \leq i \leq s$.
\end{corollary}

\noindent {\bf Proof.} It follows by using a method similar to that used in Theorem 4.1.\hfill$\blacksquare$

\noindent\begin{theorem}\label{th:ch-5} Let $k$ be an integer with $k \geq 2$ and $G\in\mathcal{G}$ be a simple graph with minimum degree $\delta \geq 2k $, if $q_{3}(G)< 4\delta-2\Delta-\frac{2(3k-1)}{\delta+1}$. Then $\tau(G)\geq k$.
\end{theorem}

\noindent {\bf Proof.} We shall prove by contradiction and assume that $\tau(G)\leq k-1.$  By (2) with $k \geq 2$, we have $\sum\limits^{s}_{i=1}r_{i} =2 \sum _{1 \leq i<j \leq s}e(V_{i},V_{j}) \leq 2ks-2(k+1).$ Let $x_{l}$ denote the multiplicity of $l$ in $\{r_{1},r_{2},\cdots,r_{s}\}$ for $1 \leq l \leq 2k-1$. By Corollary 4.3, $k \leq r_{1} \leq r_{2} \leq \ldots \leq r_{s}.$ Thus $x_{j}=0$ for $j=1,2,\ldots,k-1$. By (2), we have
$$kx_{k}+(k+1)x_{k+1}+\cdots+(2k-1)x_{2k-1}+2k(s-(x_{k}+x_{k+1}+\cdots+x_{2k-1})) \leq \sum\limits^{s}_{i=1}r_{i} \leq 2ks-2(k+1),$$
which implies that
$kx_{k}+(k-1)x_{k+1}+\cdots+2x_{2k-2}+x_{2k-1} \geq 2(k+1).$
Let $h$ be the smallest index such that $x_{h} \neq 0$, then we have
$$ (2k-h)x_{h}+(2k-h-1)x_{h+1}+\cdots+2x_{2k-2}+x_{2k-1} \geq 2(k+1).\eqno{(5)}$$

{\bf Case 1.} $x_{h} \geq 2$.

Since $h \geq k$, we have $ 2(k+1)>2(2k-h)$. Then there exists an integer $b \geq 3$ such that $(b-1)(2k-h)<2(k+1) \leq b(2k-h)$. By $2(k+1) \leq b(2k-h)$, we have $h \leq \frac{(2b-2)k-2}{b}.$ It follows by $(b-1)(2k-h)<2(k+1)$ and by (5) that $x_{h}+x_{h+1}+\cdots+x_{2k-2}+x_{2k-1} \geq b,$ and so by (1), we have $2k-1 \geq r_{b} \geq \cdots \geq r_{2} \geq r_{1}.$ By Lemma 2.6, $\mid V_{i} \mid \geq \delta+1$ with $1 \leq i \leq b.$ Let $V'=V \backslash (V_{1} \cup V_{2})$, then $\mid V' \mid \geq \mid V_{3}\mid + \cdots+ \mid V_{b} \mid \geq (b-2)(\delta+1).$ Let $y=e(V_{1},V_{2})\geq0.$ The signless Laplacian quotient matrix of $G$ with respect to the partition $(V_{p},V_{q},V')$ is
$$Q_{3}=
  \left(\begin{array}{ccc}2\overline{d_{1}}-\frac{h}{\mid V_{1}\mid} & \frac{y}{\mid V_{1}\mid} & \frac{h-y}{\mid V_{1}\mid}
\\
   \frac{y}{\mid V_{2}\mid} & 2\overline{d_{2}}-\frac{h}{\mid V_{2}\mid} & \frac{h-y}{\mid V_{2}\mid} \\
\frac{h-y}{\mid V'\mid}&\frac{h-y}{\mid V'\mid}& 2\overline{d'}-\frac{2(h-y)}{\mid V'\mid}
  \end{array}\right),
$$
where $\overline{d'}$ denotes the average degree of $V'$ in $G$, and $\overline{d_{i}}$ denotes the average degree of $V_{i}$ in $G$ for $i=1,2$.
By Lemma 2.4, we have
$q_{1}(Q_{3}) \leq 2\max\{\overline{d_{1}},\overline{d_{2}},\overline{d'}\}.$ By $tr(Q_{3})=q_{1}(Q_{3})+q_{2}(Q_{3})+q_{3}(Q_{3})$ and Lemma 2.11, we have
\begin{align*}
&q_{2}(Q_{3})+q_{3}(Q_{3})=tr(Q_{3})-q_{1}(Q_{3})\geq 2 \overline{d_{1}}+
2\overline{d_{2}}+2\overline{d'}-\frac{h}{\mid V_{1}\mid}-\frac{h}{\mid V_{2}\mid}-
\frac{2(h-y)}{\mid V'\mid} \\&- 2\max\{\overline{d_{1}},\overline{d_{2}},\overline{d'}\}
\geq 4\delta -\frac{{2\frac{b-1}{b-2}h}}{\delta+1} \geq 2( 2\delta-\frac{\frac{2 (b-1)^{2}}{b(b-2)}k-\frac{2(b-1)}{b(b-2)}}{\delta+1})>2(2\delta-\frac{3k-1}{\delta+1}).
\end{align*}
By Lemma 2.5, we have $q_{3}(Q(G)) \geq q_{3}(Q_{3}).$
Therefore, we have
\begin{align*}
q_{3}(Q(G))\geq q_{3}(Q_{3})&>2(2\delta-\frac{3k-1}{\delta+1})-
q_{2}(Q_{3}) \geq 4\delta-\frac{2(3k-1)}{\delta+1}- 2\max\{\overline{d_{1}},\overline{d_{2}},\overline{d'}\}\\&>4\delta-\frac{2(3k-1)}{\delta+1}-2\Delta,
\end{align*}
which is contrary to the assumption in Theorem 4.5 that $q_{3}(Q(G))<4\delta-2\Delta-\frac{2(3k-1)}{\delta+1}$.

{\bf Case 2.} $x_{h} = 1$.

In this case, (5) becomes $(2k-h-1)x_{h+1}+\cdots+2x_{2k-2}+x_{2k-1} \geq 2(k+1)-(2k-h)=h+2 \geq k+2.$ Let $h'$ be the smallest index such that $x_{h'}>0$ with $h'>h$, then
$$(2k-h')x_{h'}+(2k-h'-1)x_{h'+1}+\cdots+2x_{2k-2}+x_{2k-1} \geq h+2 \geq k+2.\eqno{(6)}$$
Since $ h' > h \geq k$, we have $h'+2>k$. So $k+2>2k-h'$. Thus there exists an integer $b' \geq 2$ such that $(b'-1)(2k-h')<k+2 \leq b'(2k-h')$. By $k+2 \leq b'(2k-h')$, we have $ h'\leq \frac{(2b'-1)k-2}{b'}$. It follows by $(b'-1)(2k-h')<k+2$ and by (6) that $x_{h'}+x_{h'+1}+\cdots+x_{2k-2}+x_{2k-1} \geq b',$ and so by (1), we have $2k-1 \geq r_{b'+1} \geq \cdots \geq r_{2} \geq r_{1}. $ By Lemma 2.6, $\mid V_{i} \mid \geq \delta+1$ with $1 \leq i \leq b'+1.$ Let $V'=V \backslash (V_{1} \cup V_{2})$, then $\mid V' \mid \geq \mid V_{3}\mid + \cdots+ \mid V_{b'+1} \mid \geq (b'-1)(\delta+1).$ Let $y=e(V_{1},V_{2})\geq0.$ The signless Laplacian quotient matrix of $G$ with respect to the partition $(V_{p},V_{q},V')$ is
$$Q_{3}=
  \left(\begin{array}{ccc}2\overline{d_{1}}-\frac{h}{\mid V_{1}\mid} & \frac{y}{\mid V_{1}\mid} & \frac{h-y}{\mid V_{1}\mid}
\\
   \frac{y}{\mid V_{2}\mid} & 2\overline{d_{2}}-\frac{h'}{\mid V_{2}\mid} & \frac{h'-y}{\mid V_{2}\mid} \\
\frac{h-y}{\mid V'\mid}&\frac{h'-y}{\mid V'\mid}& 2\overline{d'}-\frac{h'+h-2y}{\mid V'\mid}
  \end{array}\right),
$$
where $\overline{d'}$ denotes the average degree of $V'$ in $G$, and $\overline{d_{i}}$ denotes the average degree of $V_{i}$ in $G$ for $i=1,2$.
By Lemma 2.4, we have $q_{1}(Q_{3}) \leq 2\max\{\overline{d_{1}},\overline{d_{2}},\overline{d'}\}.$
By $tr(Q_{3})=q_{1}(Q_{3})+q_{2}(Q_{3})+q_{3}(Q_{3})$ and Lemma 2.12, we have
\begin{align*}
q_{2}(Q_{3})+q_{3}(Q_{3})&=tr(A_{3})-q_{1}(Q_{3})\geq 2\overline{d_{1}}+
2\overline{d_{2}}+2\overline{d'}-\frac{h}{\mid V_{1}\mid}-\frac{h'}{\mid V_{2}\mid}-\frac{h'+h-2y}{\mid V'\mid}\\&-2\max\{\overline{d_{1}},\overline{d_{2}},\overline{d'}\}\geq 4\delta -\frac{b'h+b'h'}{(\delta+1)(b'-1)} \geq 4\delta-\frac{2b'h'}{(\delta+1)(b'-1)}\\&\geq 2(2\delta-\frac{\frac{2b'-1}{b'-1}k-\frac{2}{b'-1}}{\delta+1}) >2(2\delta-\frac{3k-1}{\delta+1}).
\end{align*}
By Lemma 2.5, we have $q_{3}(Q(G)) \geq q_{3}(Q_{3})$.
Therefore, we have
\begin{align*}
q_{3}(Q(G))\geq q_{3}(Q_{3})&\geq 2(2\delta-\frac{3k-1}{\delta+1})-
q_{2}(Q_{3}) \geq 4\delta-\frac{2(3k-1)}{\delta+1}- 2\max\{\overline{d_{1}},\overline{d_{2}},\overline{d'}\}\\&>4\delta-\frac{2(3k-1)}{\delta+1}-2\Delta,
\end{align*}
which is contrary to the assumption in Theorem 4.5 that $q_{3}(Q(G))<4\delta-2\Delta-\frac{2(3k-1)}{\delta+1}$. This completes the proof.\hfill$\blacksquare$

\noindent\begin{corollary}\label{co:ch-6}
For $k \geq 2, $ and $d \geq 2k$. If $G\in\mathcal{G}$ is a d-regular simple graph with $q_{3}(G)< 2d-\frac{2(3k-1)}{d+1}$, then $\tau(G)\geq k$.
\end{corollary}

\noindent\begin{theorem}\label{th:ch-7} (\cite{GLLY})
For $k \geq 2,$ $ \delta \geq 2k,$ if $\lambda_{2}(G)< \delta-\frac{2k-1}{\delta+1}$. Then there exist no indices $p$ and $q$ with $1 \leq p \neq q \leq s$ such that $e(V_{p},V_{q})=0$ and $r_{p}, r_{q} \leq 2k-1$.
\end{theorem}

\noindent\begin{corollary}\label{co:ch-8}
For $k \geq 2,$ $ \delta \geq 2k,$ if $q_{2}(G)< 2\delta-\frac{2k-1}{\delta+1}$. Then there exist no indices $p$ and $q$ with $1 \leq p \neq q \leq s$ such that $e(V_{p},V_{q})=0$ and $r_{p}, r_{q} \leq 2k-1$.
\end{corollary}

\noindent {\bf Proof.} It follows by Lemma 2.9 and Theorem 4.7.\hfill$\blacksquare$

\noindent\begin{theorem}\label{th:ch-9} (\cite{GLLY})
For $k \geq 2$, $\delta \geq 2k$, if $\lambda_{2}(G)< \delta-\frac{2k-1}{\delta+1}$. Then $r_{i} \geq k$ for any $i$ with $1 \leq i \leq s$.
\end{theorem}

\noindent\begin{corollary}\label{co:ch-10}
For $k \geq 2, $ if $\delta \geq 2k$ and $q_{2}(G)< 2\delta-\frac{2k-1}{\delta+1}$. Then $r_{i} \geq k$ for any $i$ with $1 \leq i \leq s$.
\end{corollary}

\noindent {\bf Proof.} It follows by Lemma 2.9 and Theorem 4.9.\hfill$\blacksquare$

 The following Theorem 4.11 was given by Gu et al. \cite{GLLY}. In the following, we give a different proof by using the signless Laplacian quotient matrix.

\noindent\begin{theorem}\label{th:ch-11} (\cite{GLLY})
Let $k$ be an integer with $k \geq 2$ and $G$ be a simple graph with minimum degree $\delta \geq 2k $. If $q_{2}(G)< 2\delta-\frac{3k-1}{\delta+1}$, then $\tau(G)\geq k$.
\end{theorem}

\noindent {\bf Proof.} Use the similar arguments as that used in Theorem 4.5. By Corollaries 4.8 and 4.10, we have the following two cases.

{\bf Case 1.} $x_{h} \geq 2$.

In this case, we have
$$2q_{2}(Q(G))\geq q_{2}(Q_{3})+q_{3}(Q_{3})>4\delta-\frac{2(3k-1)}{\delta+1},$$
which is contrary to the assumption in Theorem 4.11 that $q_{2}(Q(G))<2\delta-\frac{3k-1}{\delta+1}$.

{\bf Case 2.} $x_{h} = 1$.

In this case, we have
\begin{align*}
2q_{2}(Q(G))&\geq q_{2}(Q_{3})+q_{3}(Q_{3})\geq 2(\overline{d_{1}}+\overline{d_{2}}+\overline{d'}-\frac{3k-1}{\delta+1})-q_{1}(Q_{3})\geq 2\overline{d_{1}}+2\overline{d_{2}}+2\overline{d'}\\&-\frac{2(3k-1)}{\delta+1}- 2\max\{\overline{d_{1}},\overline{d_{2}},\overline{d'}\}>4\delta-\frac{2(3k-1)}{\delta+1},
\end{align*}
which is contrary to the assumption in Theorem 4.11 that $q_{2}(Q(G))<2\delta-\frac{3k-1}{\delta+1}$.
This completes the proof.\hfill$\blacksquare$

By Theorem 4.11, we have the following corollary.

\noindent\begin{corollary}\label{co:ch-12}
For $k \geq 2, $ and $d \geq 2k$. If $G$ is a d-regular simple graph with $q_{2}(G)< 2d-\frac{3k-1}{d+1}$, then $\tau(G)\geq k$.
\end{corollary}

\section{Other eigenvalue conditions}
\label{sec:ch-fudi}

In this section, we will investigate the relationship between other eigenvalues and $\tau(G)$, $\kappa'(G)$ in a simple graph $G\in\mathcal{G}$. In addition, we know that Gu et al. in \cite{GLLY} gave the relationship between $\mu_{n-1}$, $\lambda_{2}$ and $\tau(G)$, $\kappa'(G)$ of a simple graph $G$, respectively.

\noindent\begin{theorem}\label{th:ch-1} Let $k$ be an integer with $k \geq 2$ , $G\in\mathcal{G}$ be a simple graph with minimum degree $\delta \geq 2k-1$.

$(i)$  If $\mu_{n-2}(G)>2\Delta-2\delta+\frac{4(k-1)}{\delta+1}$, then $\kappa'(G)\geq k$.

$(ii)$ If $q_{3}(G)<3\delta-\Delta-\frac{4(k-1)}{\delta+1}$, then $\kappa'(G)\geq k$.
\end{theorem}

\noindent {\bf Proof.} By Theorem 2.10, we have $\mu_{n-2} \leq \Delta-\lambda_{3}$, $\delta+\lambda_{3} \leq q_{3}$. It follows by Theorem 3.1.\hfill$\blacksquare$

\noindent\begin{theorem}\label{th:ch-2} Let $k$ be an integer with $k \geq 2$ and $G\in\mathcal{G}$ be a simple graph with minimum degree $\delta\geq 2k$.

$(i)$  If $\mu_{n-2}(G)> 2\Delta- 2\delta+\frac{2(3k-1)}{\delta+1}$, then $\tau(G) \geq k$.

$(ii)$ If $q_{3}(G)< 3\delta-\Delta-\frac{2(3k-1)}{\delta+1}$, then $\tau(G) \geq k$.
\end{theorem}

\noindent {\bf Proof.} By Theorem 2.10, we have $\mu_{n-2} \leq \Delta-\lambda_{3}$, $\delta+\lambda_{3} \leq q_{3}$. It follows by Theorem 4.2.\hfill$\blacksquare$

\noindent\begin{theorem}\label{th:ch-3} Let $k$ be an integer with $k \geq 2$ , $G\in\mathcal{G}$ be a simple graph with minimum degree $\delta \geq 2k-1$.

$(i)$  If $\mu_{n-2}(G)>3\Delta-3\delta+\frac{4(k-1)}{\delta+1}$, then $\kappa'(G)\geq k$.

$(ii)$ If $\lambda_{3}(G)<3\delta-2\Delta-\frac{4(k-1)}{\delta+1}$, then $\kappa'(G)\geq k$.
\end{theorem}

\noindent {\bf Proof.} By Theorem 2.10, we have $\mu_{n-2} \leq \Delta-\lambda_{3}$, $\delta+\lambda_{3} \leq q_{3}$. It follows by Theorem 3.3.\hfill$\blacksquare$

\noindent\begin{theorem}\label{th:ch-4} Let $k$ be an integer with $k \geq 2$ and $G\in\mathcal{G}$ be a simple graph with minimum degree $\delta\geq 2k$.

$(i)$  If $\mu_{n-2}(G)> 3\Delta- 3\delta+\frac{2(3k-1)}{\delta+1}$, then $\tau(G) \geq k$.

$(ii)$ If $\lambda_{3}(G)< 3\delta-2\Delta-\frac{2(3k-1)}{\delta+1}$, then $\tau(G) \geq k$.
\end{theorem}

\noindent {\bf Proof.} By Theorem 2.10, we have $\mu_{n-2} \leq \Delta-\lambda_{3}$, $\delta+\lambda_{3} \leq q_{3}$. It follows by Theorem 4.5.\hfill$\blacksquare$

\section{Spectral conditions for edge connectivity and the spanning tree packing number in multigraphs}
\label{sec:ch-fudi}

In this section, we investigate $\tau(G)$ and $\kappa'(G)$ of a multigraph $G\in\mathcal{G}$ from (signless Laplacian) spectral perspective. And we study the third small Laplacian eigenvalue conditions on $\tau(G)$ and $\kappa'(G)$. In additional, we know that Gu in \cite{Gu} gave the relationship between $\mu_{n-1}$, $\lambda_{2}$ and $\tau(G)$, $\kappa'(G)$ of a multigraph $G$.

\noindent\begin{theorem}\label{th:ch-1}
 Let $G\in\mathcal{G}$ be a multigraph with multiplicity $m$ and minimum degree $\delta \geq 2k-1 \geq 3,$ and let $l=\max\{\lceil\frac{\delta+1}{m}\rceil,2\}.$  If $\lambda_{3}(G)<2\delta-\Delta-\frac{4(k-1)}{l}$, then $\kappa'(G)\geq k$.
\end{theorem}

\noindent {\bf Proof.} Using the similar method as that used in Theorem 3.1, it follows by Lemma 2.7.\hfill$\blacksquare$

\noindent\begin{theorem}\label{th:ch-2} Let $G\in\mathcal{G}$ be a multigraph with multiplicity $m$ and minimum degree $\delta\geq 2k$ for $k \geq 2,$ and let $l=\max\{\lceil\frac{\delta+1}{m}\rceil,2\}.$  If $\lambda_{3}(G)< 2\delta-\Delta-\frac{2(3k-1)}{l}$, then $r_{i} \geq k$ for any $i$ with $1 \leq i \leq s$.
\end{theorem}

\noindent {\bf Proof.} Using the similar method as that used in Theorem 4.1, it follows by Theorem 6.1.\hfill$\blacksquare$

\noindent\begin{theorem}\label{th:ch-3}
Let $G\in\mathcal{G}$ be a multigraph with multiplicity $m$ and minimum degree $\delta\geq 2k$ for $k \geq 2,$ and let $l=\max\{\lceil\frac{\delta+1}{m}\rceil,2\}.$  If $\lambda_{3}<2\delta-\Delta-\frac{2(3k-1)}{l}$, then $\tau(G)\geq k$.
\end{theorem}

\noindent {\bf Proof.} Using the similar method as that used in Theorem 4.2, it follows by Lemma 2.7.\hfill$\blacksquare$

\noindent\begin{corollary}\label{th:ch-4} Let $G\in\mathcal{G}$ be a multigraph with multiplicity $m$ and minimum degree $\delta\geq 2k-1$ for $k \geq 2,$ and let $l=\max\{\lceil\frac{\delta+1}{m}\rceil,2\}.$

$(i)$  If $\mu_{n-2}(G)>2\Delta-2\delta+\frac{4(k-1)}{l}$, then $\kappa'(G)\geq k$.

$(ii)$ If $q_{3}(G)<3\delta-\Delta-\frac{4(k-1)}{l}$, then $\kappa'(G)\geq k$.
\end{corollary}

\noindent {\bf Proof.} By Lemma 2.10, we have $\mu_{n-2} \leq \Delta-\lambda_{3}$, $\delta+\lambda_{3} \leq q_{3}$. It follows by Theorem 6.1.\hfill$\blacksquare$

\noindent\begin{corollary}\label{th:ch-5} Let $k$ be an integer with $k \geq 2$ and $G\in\mathcal{G}$ be a multigraph with multiplicity $m$ and minimum degree $\delta\geq 2k$, and let $l=\max\{\lceil\frac{\delta+1}{m}\rceil,2\}.$

$(i)$  If $\mu_{n-2}(G)> 2\Delta- 2\delta+\frac{2(3k-1)}{l}$, then $\tau(G) \geq k$.

$(ii)$ If $q_{3}(G)< 3\delta-\Delta-\frac{2(3k-1)}{l}$, then $\tau(G) \geq k$.
\end{corollary}

\noindent {\bf Proof.} By Lemma 2.10, we have $\mu_{n-2} \leq \Delta-\lambda_{3}$, $\delta+\lambda_{3} \leq q_{3}$. It follows by Theorem 6.3.\hfill$\blacksquare$

\noindent\begin{theorem}\label{th:ch-6}
 Let $G\in\mathcal{G}$ be a multigraph with multiplicity $m$ and minimum degree $\delta \geq 2k-1 \geq 3,$ and let $l=\max\{\lceil\frac{\delta+1}{m}\rceil,2\}.$  If $q_{3}(G)<4\delta-2\Delta-\frac{4(k-1)}{l}$, then $\kappa'(G)\geq k$.
\end{theorem}

\noindent {\bf Proof.} Using the similar method as that used in Theorem 3.3, it follows by Lemma 2.7.\hfill$\blacksquare$

\noindent\begin{theorem}\label{th:ch-7} Let $G\in\mathcal{G}$ be a multigraph with multiplicity $m$ and minimum degree $\delta\geq 2k$ for $k \geq 2,$ and let $l=\max\{\lceil\frac{\delta+1}{m}\rceil,2\}.$  If $q_{3}<4\delta-2\Delta-\frac{2(3k-1)}{l}$, then $r_{i} \geq k$ for any $i$ with $1 \leq i \leq s$.
\end{theorem}

\noindent {\bf Proof.} Using the similar method as that used in Theorem 4.1, it follows by Theorem 6.6.\hfill$\blacksquare$

\noindent\begin{theorem}\label{th:ch-8}
Let $G\in\mathcal{G}$ be a multigraph with multiplicity $m$ and minimum degree $\delta\geq 2k$ for $k \geq 2,$ and let $l=\max\{\lceil\frac{\delta+1}{m}\rceil,2\}.$  If $q_{3}<4\delta-2\Delta-\frac{2(3k-1)}{l}$, then $\tau(G)\geq k$.
\end{theorem}

\noindent {\bf Proof.} Using the similar method as that used in Theorem 4.5, it follows by Lemma 2.7.\hfill$\blacksquare$

\noindent\begin{corollary}\label{th:ch-9} Let $G\in\mathcal{G}$ be a multigraph with multiplicity $m$ and minimum degree $\delta\geq 2k-1$ for $k \geq 2,$ and let $l=\max\{\lceil\frac{\delta+1}{m}\rceil,2\}.$

$(i)$  If $\mu_{n-2}(G)>3\Delta-3\delta+\frac{4(k-1)}{l}$, then $\kappa'(G)\geq k$.

$(ii)$ If $\lambda_{3}(G)<3\delta-2\Delta-\frac{4(k-1)}{l}$, then $\kappa'(G)\geq k$.
\end{corollary}

\noindent {\bf Proof.} By Lemma 2.10, we have $\mu_{n-2} \leq \Delta-\lambda_{3}$, $\delta+\lambda_{3} \leq q_{3}$. It follows by Theorem 6.6.\hfill$\blacksquare$

\noindent\begin{corollary}\label{th:ch-10} Let $k$ be an integer with $k \geq 2$ and $G\in\mathcal{G}$ be a multigraph with multiplicity $m$ and minimum degree $\delta\geq 2k,$ and let $l=\max\{\lceil\frac{\delta+1}{m}\rceil,2\}.$

$(i)$  If $\mu_{n-2}(G)> 3\Delta- 3\delta+\frac{2(3k-1)}{l}$, then $\tau(G) \geq k$.

$(ii)$ If $\lambda_{3}(G)< 3\delta-2\Delta-\frac{2(3k-1)}{l}$, then $\tau(G) \geq k$.
\end{corollary}

\noindent {\bf Proof.} By Lemma 2.10, we have $\mu_{n-2} \leq \Delta-\lambda_{3}$, $\delta+\lambda_{3} \leq q_{3}$. It follows by Theorem 6.8.\hfill$\blacksquare$

{\bf Remark 2.} Let $G_{1}$ and $G_{2}$ be two graphs shown in Figures 2 and 3. By Matlab calculation, we have $\lambda_{3}(G_{1})\approx 1.562<2\times3-3-\frac{4(2-1)}{3+1}=2,$ $\kappa'(G_{1})=2=k,$ and $\lambda_{3}(G_{2})\approx 0.618<2\times3-4-\frac{4(2-1)}{3+1}=1,$ $\kappa'(G_{2})=2=k.$  Notice that $G_{1},G_{2}$ satisfy the conditions of Theorem 3.1 but $G_{1},G_{2}\notin \mathcal{G}$.

\begin{center}
\includegraphics [width=5 cm, height=4 cm]{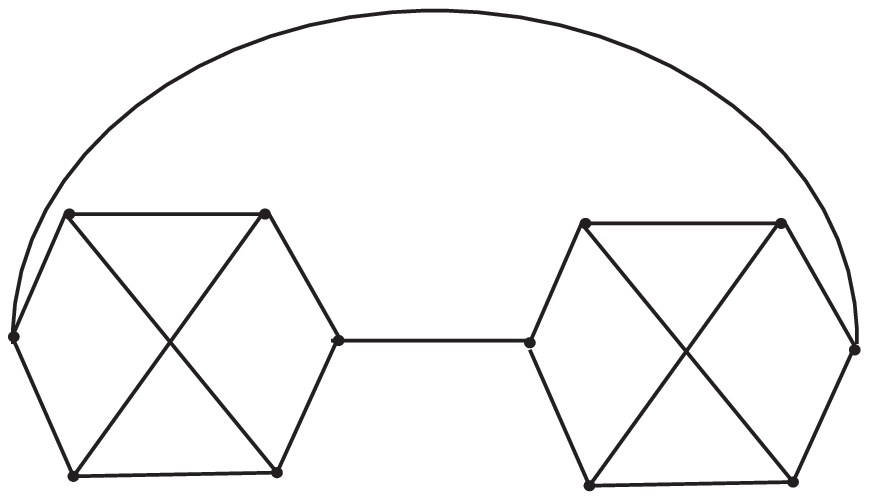}
\centerline{Figure 2: A 3-regular graph $G_{1}$ with $\lambda_{3}\approx 1.562<2.$}
\end{center}

\begin{center}
\includegraphics [width=5 cm, height=3 cm]{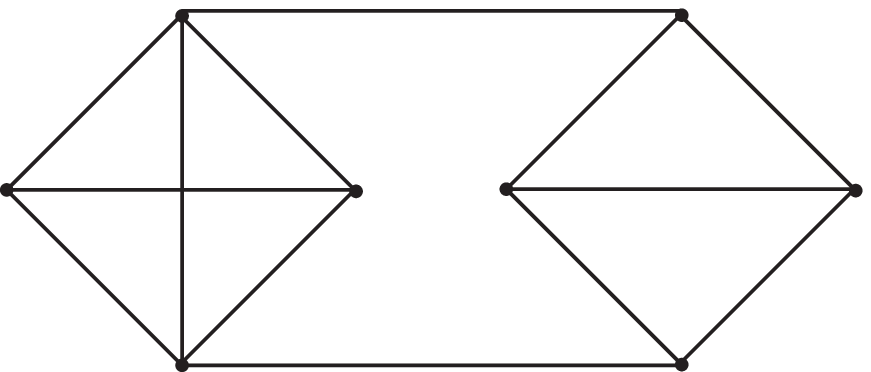}
\centerline{Figure 3: A simple graph $G_{2}$ with $\lambda_{3}\approx 0.618<1.$}
\end{center}

In this paper, we only consider graphs, which have at least two non-empty disjoint proper subsets $V_{1},V_{2}\subseteq V(G)$ such that $\kappa'(G)=e(V_{i},V\backslash V_{i})$ for $1\leq i \leq 2$ and $V\setminus(V_{1} \cup V_{2})\neq \phi.$ The proofs of Theorems 3.1 and 3.3 cannot be applied to a graph $G\notin\mathcal{G}$ since $V(G)$ in the proof may not be able to divide into a $3 \times 3$ quotient matrix which we need. We consider that the results of this paper may also be true for a graph $G\notin\mathcal{G}$, for example, see Figures 2 and 3.



\end{document}